\documentclass[11pt]{article}
\usepackage{amssymb,amsmath}
\usepackage{arrows}
\usepackage[dvips]{graphicx}
\numberwithin{equation}{section}
\newcommand{\real}{{\mathbb R}}

\renewcommand{\natural}{{\mathbb N}}

\newcommand{\Ker}{\mathrm{Ker}}
\newcommand{\Ran}{\mathrm{Ran}}

\renewcommand{\d}{\,{\rm d}}            
\newcommand{\D}{{\rm d}}                

\newcommand{\weakto}{\rightharpoonup}

\newcommand{\cL}{{\cal L}}
\newcommand{\cM}{{\cal M}}

\newcommand{\cO}{{\cal O}}

\newcommand{\cS}{{\cal S}}

\newtheorem{theorem}{Theorem}[section]
\newtheorem{lemma}[theorem]{Lemma}

\newtheorem{proposition}[theorem]{Proposition}
\newtheorem{corollary}[theorem]{Corollary}
\newtheorem{remark}[theorem]{Remark}

\newcommand{\proof}{{\noindent \bf Proof.\ }}

\newcommand{\app}{{\rm app}}

\def\build#1_#2^#3{\mathrel{
  \mathop{\kern 0pt#1}\limits_{#2}^{#3}}}

\def\QED{\mbox{}\hfill$\Box$}

\renewcommand{\:}{\thinspace :}

\newdimen\texpscorrection
\newdimen\figcenter
\def\figurewithtex #1 #2 #3 #4 #5\cr{\null
  {\goodbreak\figcenter=\hsize\relax
  \advance\figcenter by -#4truecm
  \divide\figcenter by 2
  \begin{figure}[hbt]
  \vskip #3truecm\noindent\hskip\figcenter
  \includegraphics{#1}{\hskip\texpscorrection\input #2 }
  \vskip 0.8truecm{\baselineskip=0.8\baselineskip
  \noindent \vbox{\noindent {\footnotesize #5}}\par}
  \end{figure}}}
\def\point#1 #2 #3 {\rlap{\kern #1 truecm
\raise #2 truecm \hbox{#3}}}

\begin{document}

\title{Interacting vortex pairs in inviscid and viscous planar flows}

\author{{\bf Thierry Gallay} \\[1mm] 
Universit\'e de Grenoble I\\
Institut Fourier, UMR CNRS 5582\\
B.P. 74\\
F-38402 Saint-Martin-d'H\`eres, France\\
{\tt Thierry.Gallay@ujf-grenoble.fr}}

\maketitle

\begin{abstract}
The aim of this contribution is to make a connection between two
recent results concerning the dynamics of vortices in incompressible 
planar flows. The first one is an asymptotic expansion, in the vanishing
viscosity limit, of the solution of the two-dimensional Navier-Stokes
equation with point vortices as initial data. In such a situation, it
is known \cite{Ga} that the solution behaves to leading order like a
linear superposition of Oseen vortices whose centers evolve according
to the point vortex system, but higher order corrections can also be
computed which describe the deformation of the vortex cores due to
mutual interactions. The second result is the construction by D. Smets
and J. van Schaftingen of ``desingularized'' solutions of the
two-dimensional Euler equation \cite{SVS}. These solutions are
stationary in a uniformly rotating or translating frame, and converge
either to a single vortex or to a vortex pair as the size parameter
$\epsilon$ goes to zero. We consider here the particular case of a 
pair of identical vortices, and we show that the solution of the weakly 
viscous Navier-Stokes equation is accurately described at time $t$
by an approximate steady state of the rotating Euler equation which
is a desingularized solution in the sense of \cite{SVS} with Gaussian
profile and size $\epsilon = \sqrt{\nu t}$.
\end{abstract}

\section{Introduction} \label{intro} 
Numerical simulations of freely decaying turbulence show that vortex
interactions play a crucial role in the dynamics of two-dimensional
viscous flows \cite{MW84,MW90}. In particular, vortex mergers are
responsible for the appearance of larger and larger structures in the
flow, a process which is directly related to the celebrated ``inverse
energy cascade'' \cite{Fr95}.  Although nonperturbative interactions
such as vortex mergers are extremely complex and desperately hard to
analyze from a mathematical point of view \cite{LD02,Meu05}, rigorous
results can be obtained in the perturbative regime where the distances
between the vortex centers are large compared to the typical size of
the vortex cores.

As an example of this situation, consider the case where the
initial flow is a superposition of $N$ point vortices. This means that
the initial vorticity $\omega_0$ satisfies
\begin{equation}\label{omega0}
  \omega_0 \,=\, \sum_{i=1}^N \alpha_i \,\delta(\cdot - x_i)~,
\end{equation}
where $x_1,\dots,x_N \in \real^2$ are the initial positions and
$\alpha_1,\dots,\alpha_N \in \real$ the circulations of the 
vortices. Let $\omega(x,t)$ be the solution of the two-dimensional
vorticity equation
\begin{equation}\label{Veq}
  \frac{\partial\omega}{\partial t} + u \cdot \nabla \omega \,=\, 
  \nu \Delta \omega~, \qquad x \in \real^2~, \quad t > 0~, 
\end{equation}
with initial data $\omega_0$, where $u(x,t)$ is the velocity 
field defined by the Biot-Savart law
\begin{equation}\label{BS}
  u(x,t) \,=\, \frac{1}{2\pi} \int_{\real^2}
  \frac{(x - y)^{\perp}}{|x -y|^2} \,\omega(y,t)\d y~,
  \qquad x \in \real^2~, \quad t > 0~.
\end{equation}
Solutions of \eqref{Veq}, \eqref{BS} with singular initial data of the
form \eqref{omega0} were first constructed by Benfatto, Esposito, and
Pulvirenti \cite{BEP85}. More generally, if $\omega_0 \in
\cM(\real^2)$ is any finite measure, Giga, Miyakawa, and Osada
\cite{GMO88} have shown that the vorticity equation \eqref{Veq} has a
global solution with initial data $\omega_0$, which moreover is unique
if the total variation norm of atomic part of $\omega_0$ is small
compared to the kinematic viscosity $\nu$. This last restriction has
been removed recently by I.~Gallagher and the author \cite{GG05}, so
we know that \eqref{Veq} has a unique global solution $\omega \in
C^0((0,\infty),L^1(\real^2) \cap L^\infty(\real^2))$ with initial data
\eqref{omega0}, no matter how small the viscosity coefficient is. This
solution is uniformly bounded in $L^1(\real^2)$, and the total
circulation $\int_{\real^2}\omega(x,t)\d x$ is conserved.

In the vanishing viscosity limit, the motion of point vortices
in the plane is described by a system of ordinary differential 
equations introduced by Helmholtz \cite{He58} and Kirchhoff 
\cite{Ki76}. If $z_1(t),\dots,z_N(t) \in \real^2$ denote the 
positions of the vortices, the system reads
\begin{equation}\label{PW}
  \frac{\D}{\D t}z_i(t) \,=\, \frac{1}{2\pi} \sum_{j\neq i} \alpha_j 
  \,\frac{(z_i(t) - z_j(t))^\perp}{|z_i(t)-z_j(t)|^2}~, 
  \qquad i = 1,\dots,N~,
\end{equation}
and the initial conditions $z_i(0) = x_i$ for $i = 1,\dots,N$ are
determined by \eqref{omega0}. A lot is known about the dynamics 
of the {\em point vortex system} \eqref{PW}, see e.g. \cite{New01}
for a recent monograph devoted to this problem. Most remarkably, 
\eqref{PW} is a Hamiltonian system with $N$ degrees of freedhom, 
which always possesses three independent involutive first 
integrals. In particular, system \eqref{PW} is {\em integrable} 
if $N \le 3$, whatever the vortex circulations $\alpha_1,\dots,\alpha_N$
may be. In the simple situation where $N = 2$, both vortices 
rotate with constant angular speeed around the common vorticity 
center, see Fig.\thinspace 1\thinspace(left). In the exceptional case 
where $\alpha_1 + \alpha_2 = 0$, there is no center of vorticity and 
the vortices move with constant speed along parallel straight lines, see 
Fig.\thinspace 1\thinspace(right). 
\begin{figure}
\begin{tabular}{cc}
\hspace{1cm}\includegraphics[width=6.6cm,height=6.0cm]{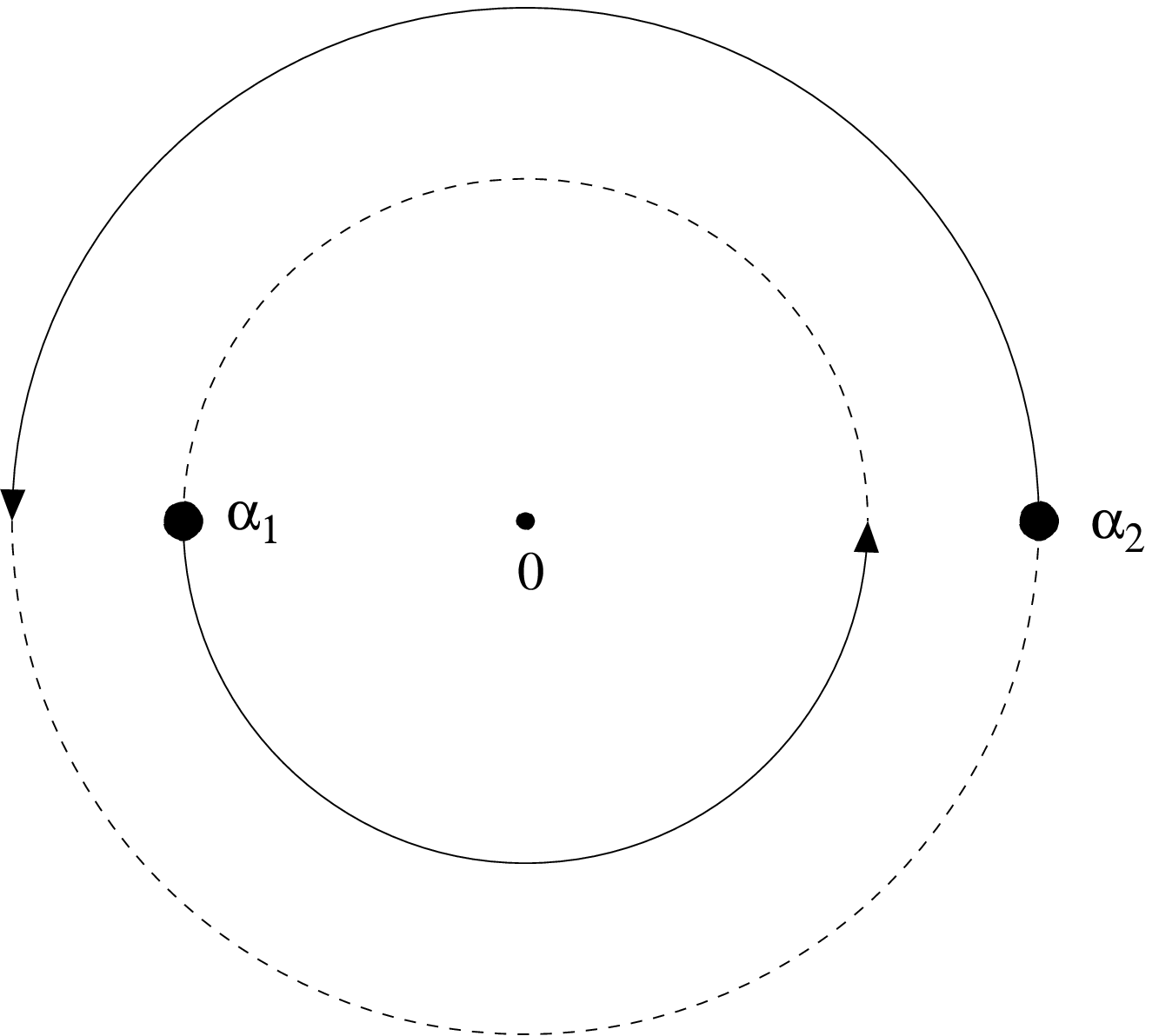} & 
\hspace{2cm}\includegraphics[width=4.2cm,height=4.5cm]{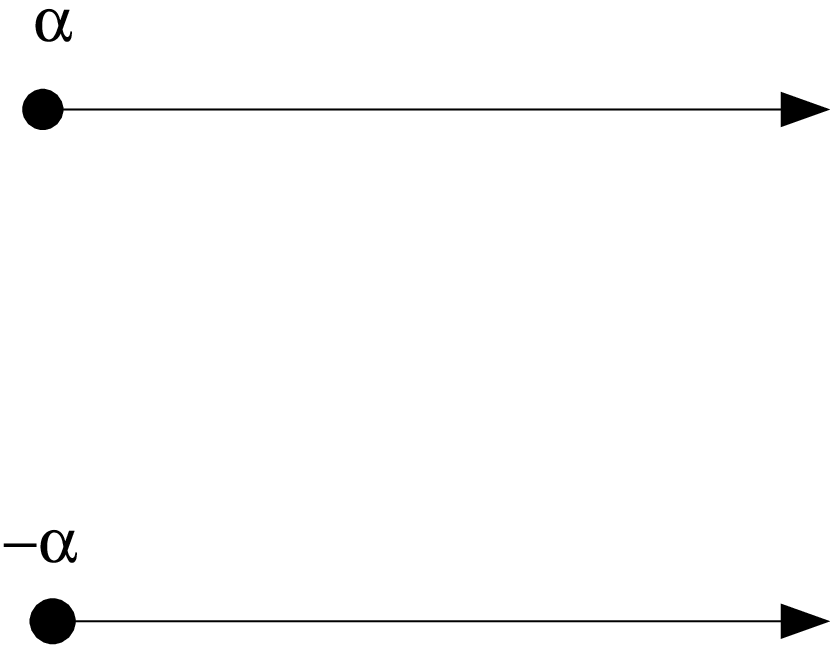} \\
\end{tabular}
\caption{\small The motion of two point vortices whith circulations
$\alpha_1 > \alpha_2 > 0$ (left) and $\alpha_1 + \alpha_2 = 0$ (right).}
\label{Fig1}
\end{figure}

It should be remarked that system \eqref{PW} is not always 
globally well-posed: if $N \ge 3$ and if the circulations 
$\alpha_1,\dots,\alpha_N$ are not all of the same sign, vortex
collisions may occur in finite time for exceptional initial 
configurations \cite{MP94,New01}. In what follows, we always assume
that system \eqref{PW} is well-posed on some time interval 
$[0,T]$, and we denote by $d$ the {\em minimal distance} between 
any two vortices on this interval:
\begin{equation}\label{dmin}
  d \,=\, \min_{t \in [0,T]}\,\min_{i\neq j}\,|z_i(t) - z_j(t)|
  \,>\, 0~. 
\end{equation}
We also introduce the {\em turnover time} 
\begin{equation}\label{T0def}
  T_0 \,=\, \frac{d^2}{|\alpha|}~, \quad \hbox{where} \quad
  |\alpha| \,=\, |\alpha_1| + \dots + |\alpha_N|~,
\end{equation}
which is a natural time scale for the inviscid dynamics described 
by \eqref{PW}. For instance, for a pair of vortices with circulations
of the same sign, one can check that the rotation period of each
vortex around the center is $4\pi^2 T_0$. 

When the viscosity $\nu$ is nonzero, the point vortices in the initial 
data \eqref{omega0} are smoothed out by diffusion, and the solution 
$\omega(x,t)$ of \eqref{Veq} is no longer described by the point
vortex system. In the particular case where $N = 1$, the unique 
solution of \eqref{Veq} with initial data $\omega_0 = \alpha 
\delta_0$ is the {\em Lamb-Oseen vortex:}
\begin{equation}\label{Oseen}
  \omega(x,t) \,=\, \frac{\alpha}{\nu t}\,G\Bigl(\frac{x}
  {\sqrt{\nu t}}\Bigr)~, \qquad u(x,t) \,=\,\frac{\alpha}{\sqrt{\nu t}}
  \,v^G\Bigl(\frac{x}{\sqrt{\nu t}}\Bigr)~,
\end{equation}
where the vorticity profile $G$ and the velocity profile $v^G$ 
have the following explicit expressions:
\begin{equation}\label{Gdef}
  G(\xi) \,=\, \frac{1}{4\pi} \,e^{-|\xi|^2/4}~,\qquad
  v^G(\xi) \,=\, \frac{1}{2\pi}\frac{\xi^\perp}{|\xi|^2}
  \Bigl(1 -  e^{-|\xi|^2/4}\Bigr)~.
\end{equation}
As was shown by C.E.~Wayne and the author, Oseen vortices describe the 
long-time asymptotics of all solutions of the two-dimensional 
Navier-Stokes equation for which the vorticity distribution is 
integrable, and are also the only self-similar solutions of this
equation with integrable vorticity profile \cite{GW05}. 

When $N \ge 2$, the solution of \eqref{Veq} with initial data 
\eqref{omega0} is not explicit, but in some parameter 
regimes it can be approximated by a linear superposition of Oseen 
vortices whose centers evolve according to the point vortex 
system \eqref{PW}. More precisely, we have the following result: 

\begin{theorem}\label{thm1} {\bf \cite{Ga}} Given pairwise distinct
initial positions $x_1,\dots,x_N \in \real^2$ and nonzero circulations 
$\alpha_1,\dots,\alpha_N \in \real$, fix $T > 0$ such that
the point vortex system \eqref{PW} is well-posed on the time 
interval $[0,T]$. Let $d > 0$ be the minimal distance \eqref{dmin}
and $T_0$ the turnover time \eqref{T0def}. Then the (unique) 
solution of the two-dimensional vorticity equation \eqref{Veq} 
with initial data \eqref{omega0} satisfies
\begin{equation}\label{thm1app}
  \frac{1}{|\alpha|}\int_{\real^2}\Bigl|\omega(x,t) - \sum_{i=1}^N 
  \frac{\alpha_i}{\nu t}\,G\Bigl(\frac{x-z_i(t)}{\sqrt{\nu t}}\Bigr)
  \Bigr|\d x \,\le \, K\,\frac{\nu t}{d^2}~, \qquad t \in (0,T]~,
\end{equation}
where $z(t) = \{z_1(t),\dots,z_N(t)\}$ is the solution of \eqref{PW} 
and $K$ is a (dimensionless) constant depending only on 
the ratio $T/T_0$. 
\end{theorem}

Theorem~\ref{thm1} gives nontrivial information on the solution of
\eqref{Veq} in the {\em weak interaction regime}, where the size
$\cO(\sqrt{\nu t})$ of the vortex cores is much smaller than the
distance $d$ between the centers. If the initial data are fixed, a
convenient way to achieve this is to assume that the viscosity $\nu$
is small compared to the total circulation $|\alpha|$, and that the
observation time $T$ is smaller than or comparable to the turnover
time $T_0$. Indeed, using definition \eqref{T0def}, we see that
\[
  \frac{\nu t}{d^2} \,=\, \frac{\nu}{|\alpha|}\,\frac{t}{T_0}~.
\]
In particular, in the vanishing viscosity limit, we obtain 
the following corollary.

\begin{corollary}\label{cor1} {\bf \cite{Ga}}
Under the assumptions of Theorem~\ref{thm1}, the solution 
$\omega^\nu(x,t)$ of the viscous vorticity equation \eqref{Veq} 
with initial data \eqref{omega0} satisfies
\begin{equation}\label{weakconv}
  \omega^\nu(\cdot,t) 
  ~\xrightharpoonup[\nu \to 0]{\hbox to 8mm{}}~ 
  \sum_{i=1}^N \alpha_i \,\delta(\cdot - z_i(t))~, \quad 
  \hbox{for all } t \in [0,T]~,
\end{equation}
where $z(t) = \{z_1(t),\dots,\,z_N(t)\}$ is the solution of \eqref{PW}. 
\end{corollary}

In other words, the solution of the vorticity equation \eqref{Veq}
with initial data \eqref{omega0} converges weakly, in the vanishing
viscosity limit, to a superposition of point vortices which evolve
according to the point vortex dynamics \eqref{PW}. In particular, 
Corollary~\ref{cor1} provides a natural and rigorous derivation 
of the point vortex system itself from the Navier-Stokes equation. 
As is well-known, system \eqref{PW} can also be rigorously derived
from Euler's equation \cite{Ma88,MP93}, but in the latter approach it 
is necessary to regularize the initial vorticity because we do not
know how to solve Euler's equation with singular data such as 
\eqref{omega0}. In this respect, it is important to keep in mind that 
the right-hand side of \eqref{weakconv} is {\em not} a weak solution of 
the inviscid vorticity equation $\partial_t \omega + u \cdot \nabla
\omega = 0$. For completeness, we also mention that the vanishing 
viscosity limit for solutions of \eqref{Veq} with concentrated vorticity 
has been studied by Marchioro \cite{Ma98}, who obtained the analog of 
Corollary~\ref{cor1} in that context. 

The conclusion of Theorem~\ref{thm1} can be interpreted in the
following simple and somewhat naive way: If we solve the
two-dimensional vorticity equation \eqref{Veq} with point vortices as
initial data, the diffusion term $\nu \Delta\omega$ in \eqref{Veq}
smooths out the point vortices into Oseen vortices, and the advection
term $u\cdot\nabla\omega$ translates the vortex centers according to
the point vortex dynamics \eqref{PW}.  While correct, this
interpretation ignores the important fact that the advection of a
smooth vortex by the inhomogeneous velocity field created by the other
vortices results non only in a translation of the vortex center, but
also in a {\em deformation} of the vortex core. In our case, this
deformation is of order $\cO(\nu t/d^2)$ in $L^1$ norm, and therefore
does not appear in \eqref{thm1app} because it is included in the
error term. However, as is explained in \cite{Ga}, such a small
deformation of the vortex profile creates a {\em self-interaction}
effect of size $\cO(1)$, which basically counterbalances the influence
of the external velocity field, except for a rigid translation. 

As a matter of fact, self-interactions play a crucial role in the
proof of Theorem~\ref{thm1}, and we even believe that is it not
possible to establish \eqref{thm1app} without computing a higher order
approximation of the solution. A systematic asymptotic expansion is 
carried out in \cite{Ga}, but if we only keep the leading order nonradial
corrections to the Oseen vortices we obtain the following approximate
solution of \eqref{Veq}:
\begin{equation}\label{om_app}
  \omega_\app(x,t) \,=\, \sum_{i=1}^N \frac{\alpha_i}{\nu t}
  \left\{G\Bigl(\frac{x-z_i(t)}{\sqrt{\nu t}}\Bigr) + 
  \frac{\nu t}{d^2}\,F_i\Bigl(\frac{x-z_i(t)}{\sqrt{\nu t}}\,,t\Bigr) 
  \right\}~,
\end{equation} 
where
\begin{equation}\label{Fidef}
  F_i(\xi,t) \,=\, a(|\xi|) \sum_{j\neq i} \frac{\alpha_j}{\alpha_i}
  \,\frac{d^2}{|z_i(t)-z_j(t)|^2}
  \Bigl(2\frac{|\xi \cdot (z_i(t)-z_j(t))|^2}{|\xi|^2 |z_i(t) - 
  z_j(t)|^2} - 1\Bigr) + \cO\Bigl(\frac{\nu}{|\alpha|}\Bigr)~. 
\end{equation}
Here $a : (0,\infty) \to \real$ is a smooth, positive function 
satisfying $a(r) \approx C_1 r^2$ as $r \to 0$ and $a(r) \approx 
C_2 r^4 e^{-r^2/4}$ as $r \to \infty$ for some $C_1, C_2 > 0$. 
In the regime where the viscosity $\nu$ is much smaller than the
circulations $\alpha_i$ of the vortices, we can use formulas 
\eqref{om_app}, \eqref{Fidef} to compute, to leading order in 
our expansion parameter $\nu t/d^2$, the deformation of the vortex
cores due to mutual interactions. For any fixed $t \in (0,T]$, 
this first order correction depends only on the relative positions
of the vortex centers, which are determined by \eqref{PW}. In the 
particular case of a pair of vortices with equal or opposite 
circulations, the level lines of the vorticity distribution 
$\omega_\app(x,t)$ are represented in Fig.~2. 
\begin{figure}
\begin{tabular}{cc}
\hspace{1.0cm}\includegraphics[width=6.6cm,height=5.0cm]{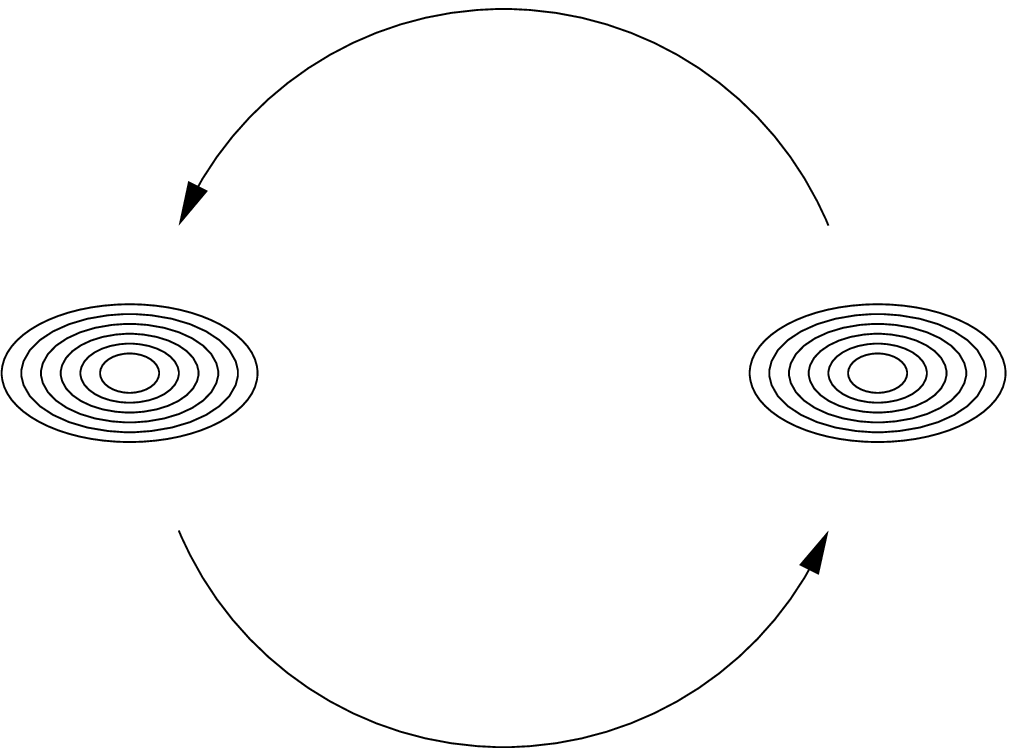} & 
\hspace{2.0cm}\includegraphics[width=4.5cm,height=5.0cm]{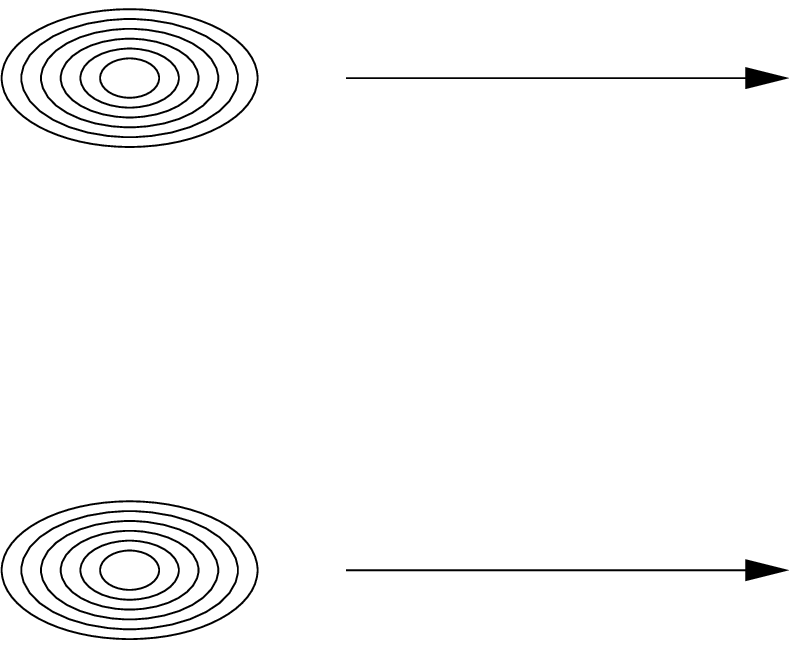} \\
\end{tabular}
\caption{\small The level lines of the vorticity distribution for 
a pair of vortices with equal (left) or opposite (right) circulations.}
\label{Fig2}
\end{figure}

The discussion above shows that the dynamics of weakly interacting
viscous vortices is essentially driven by two different mechanisms:
{\em diffusion}, which is responsible for the continuous growth of the
vortex cores, and {\em advection}, which creates the motion of the
vortex centers and the deformation of the vortex profiles.  The latter
effect persists in the vanishing viscosity limit, and it is therefore
reasonable to expect that, if we can find solutions of Euler's
equation describing widely separated Gaussian vortices, these inviscid
solutions will provide an accurate approximation of the viscous
$N$-vortex solution considered in Theorem~\ref{thm1}, if $\nu$ is
sufficiently small. The aim of this contribution is to explore this
idea in the particular case of a single {\em vortex pair}. In this
simple situation, we only need to consider solutions of Euler's
equation that are stationary in a uniformly rotating or translating
frame.

Assume thus that $N = 2$ and, for definiteness, that both circulations
$\alpha_1, \alpha_2$ are positive. Given $r_1, r_2 > 0$ such that 
$\alpha_1 r_1 = \alpha_2 r_2$, let $d = r_1 + r_2$ and $\Omega = 
(\alpha_1 + \alpha_2)/(2\pi d^2)$. We consider the inviscid vorticity 
equation in a rotating frame with angular speed $\Omega$\:
\begin{equation}\label{Vrot}
  \partial_t \omega + (u - \Omega x^\perp) \cdot \nabla 
  \omega \,=\,  0~.
\end{equation}
Formally, the vorticity distribution
\begin{equation}\label{twov}
  \omega_0 \,=\, \alpha_1 \delta(\cdot - x_1) + \alpha_2
  \delta(\cdot - x_2)~, \qquad \hbox{where} \quad x_1 = 
  \begin{pmatrix} r_1 \\ 0 \end{pmatrix}~, \quad x_2 = 
  \begin{pmatrix} -r_2 \\ 0 \end{pmatrix}~, 
\end{equation}
is a stationary solution of \eqref{Vrot}. In the laboratory 
frame, this corresponds to a time periodic solution of \eqref{PW} 
where both vortices rotate around the origin with angular velocity
$\Omega$. Now, let $w_* \in \cS(\real^2)$ be a nonnegative vorticity
profile which is radially symmetric, decreasing along rays, and 
normalized in the sense that $\int_{\real^2} w_*(x)\d x = 1$. Given 
$\epsilon > 0$, we look for a stationary solution of \eqref{Vrot} of 
the form
\begin{equation}\label{Vpair}
  \omega_\epsilon(x) \,=\, \frac{\alpha_1}{\epsilon^2} \,w_{1,\epsilon}
  \Bigl(\frac{x-x_{1,\epsilon}}{\epsilon}\Bigr) + \frac{\alpha_2}
  {\epsilon^2}\,w_{2,\epsilon}\Bigl(\frac{x-x_{2,\epsilon}}{\epsilon}
  \Bigr)~,
\end{equation}
where $x_{i,\epsilon} \to x_i$ and $w_{i,\epsilon} \to w_*$ as $\epsilon
\to 0$. Existence of such "desingularized" solutions of Euler's equation
has been investigated in a recent work by D. Smets and J. Van 
Schaftingen \cite{SVS}. In fact, the authors of \cite{SVS} do not 
consider rotating vortex pairs of the form \eqref{Vpair}, but they treat
a variety of other interesting cases, including a single stationary vortex
in a bounded or unbounded domain, a rotating vortex in a disk, and 
a translating vortex pair in the plane. They use the "stream function 
method", which consists in constructing (by variational methods) a
nontrivial solution to an elliptic equation of the form $-\Delta \psi = 
f_\epsilon(\psi + \frac12 \Omega |x|^2)$, where $f_\epsilon$ is a 
power-like nonlinearity which depends on $\epsilon$ in an appropriate
way. The vorticity $\omega = -\Delta \psi$ is then a stationary solution 
of \eqref{Vrot}. For technical reasons, the desingularized vorticity profiles 
obtained in \cite{SVS} are always compactly supported, but it is reasonable
to expect that similar results can be obtained with Gaussian profiles too.
We hope to clarify this issue and to extend the results of \cite{SVS}
to pairs of vortices of the same sign in a future work. 

In Section~\ref{sec2} below, we study in detail the case of two
identical vortices ($\alpha_1 = \alpha_2$). For a large class of
radially symmetric profiles $w_*$, we prove the existence of
approximate stationary solutions of \eqref{Vrot} of the form
\eqref{Vpair}.  In other words, we construct an asymptotic expansion
in powers of $\epsilon$ of the vorticity distribution \eqref{Vpair} as
a steady state of \eqref{Vrot}. Under natural symmetry
assumptions, we show that this expansion can be performed to
arbitrarily high order. Then, in Section~\ref{sec3}, we prove that the
solution $\omega^\nu(x,t)$ of the rotating viscous vorticity equation
$\partial_t \omega + (u - \Omega x^\perp) \cdot \nabla \omega = \nu
\Delta\omega$ with initial data \eqref{twov} is very close to the
inviscid stationary solution $\omega_\epsilon$ with asymptotic profile
$w_* = G$, if $\epsilon = \sqrt{\nu t}$ is sufficiently small. This
means that, when the viscosity is small, the solution of \eqref{Veq}
with initial data \eqref{twov} slowly travels through a family of
uniformly rotating solutions of Euler's equation, whose vorticity
profiles are approximately Gaussian and evolve diffusively. We expect
a similar picture to be relevant in the general situation considered
in Theorem~\ref{thm1}, although the corresponding inviscid solutions
may be more difficult to identify in that case.

\medskip\noindent{\bf Acknowledgements.} I am grateful to 
Felix Otto for suggesting the point of view adopted in this 
paper, and to Didier Smets for many stimulating discussions. 

\section{Approximate steady states of Euler's equation}
\label{sec2}

The aim of this section is to construct an asymptotic expansion
for a family of stationary solutions of the inviscid vorticity
equation, which correspond to weakly interacting vortex pairs. For
simplicity, we only consider the particular case where both vortices
have the same circulation $\alpha > 0$. As is explained in the
introduction, weak interaction means that the distance $d$ between the
vortex centers is large compared to the size of the vortex cores. If
the vorticity distribution is given by \eqref{Vpair}, this condition
is satisfied if $\epsilon > 0$ is sufficiently small. Here, we find it
more convenient to fix the size of the vortex cores, and to assume
that the distance $d$ between the centers is large. This alternative
point of view is of course equivalent to the previous one, up to a
rescaling. Note, however, that the rotation speed $\Omega$ will now
depend on $d$ and behave like $\alpha/(\pi d^2)$ as $d \to \infty$.

From now on, we fix $\alpha > 0$, $d \gg 1$, and we look for a 
stationary solution $\omega$ of \eqref{Vrot} describing a pair 
of identical vortices with circulation $\alpha$. We make the 
ansatz
\begin{equation}\label{om1}
  \omega(x) \,=\, \alpha w(x-x_d) + \alpha w(-x-x_d)~, \quad
  u(x) \,=\, \alpha v(x-x_d) - \alpha v(-x-x_d)~, 
\end{equation}
where $x_d = (d/2,0)$, $w$ is a localized vorticity profile
to be determined, and $v = K[w]$ is the velocity field obtained 
from $w$ via the Biot-Savart law \eqref{BS}. We assume that $w$ 
belongs to the Schwartz class $\cS(\real^2)$, is nonnegative, 
and satisfies the normalization condition 
\begin{equation}\label{wnorm}
  \int_{\real^2} w(x)\d x \,=\, 1~.
\end{equation}
We also impose the following symmetry
\begin{equation}\label{wsym}
  w(x_1,-x_2) \,=\, w(x_1,x_2)~, \quad
  v_1(x_1,-x_2) \,=\, -v_1(x_1,x_2)~, \quad
  v_2(x_1,-x_2) \,=\, v_2(x_1,x_2)~,
\end{equation}
which implies that $\omega(-x_1,x_2) = \omega(x_1,-x_2) = \omega(x_1,x_2)$ 
for all $x = (x_1,x_2) \in \real^2$. Finally, without loss of 
generality, we assume that
\begin{equation}\label{wmom}
  \int_{\real^2} x_1 w(x) \d x \,=\, 0~.
\end{equation} 
This means that the vorticity distribution $\omega$ defined in 
\eqref{om1} is indeed a superposition of two localized vortices 
centered at the points $\pm x_d$. 

The distribution $\omega$ will be a stationary solution of the rotating 
vorticity equation \eqref{Vrot} if the profile $w$ satisfies
\[
  \Bigl(\alpha v(x-x_d) - \alpha v(-x-x_d) - \Omega x^\perp\Bigr)
  \cdot \nabla w(x-x_d) \,=\, 0~, \qquad x \in \real^2~.
\]
Replacing $x$ by $x+x_d$ and denoting $\tilde \Omega = \Omega/\alpha$, 
we obtain the equivalent equation
\begin{equation}\label{weq}
  \Bigl(v(x) - v(-x-2x_d) - \tilde\Omega (x+x_d)^\perp\Bigr)
  \cdot \nabla w(x) \,=\, 0~, \qquad x \in \real^2~.  
\end{equation}
Note that the rotating term $\tilde\Omega (x+x_d)^\perp$ behaves 
like the velocity field $v(x)$ when $x_2$ is changed into 
$-x_2$; this implies that the symmetry \eqref{wsym} is
indeed compatible with Eq.~\eqref{weq}. To determine the rotation
speed $\Omega$, we multiply both members of \eqref{weq} by $x_2$ and
we integrate by parts over $\real^2$. We easily obtain
\[
  \tilde \Omega \int_{\real^2} (x_1+d/2) w(x)\d x \,=\, 
  \int_{\real^2}\Bigl(v_2(x) - v_2(-x-2x_d)\Bigr)w(x)\d x~.
\]
This identity can be simplified if we use \eqref{wnorm}, 
\eqref{wmom}, and the fact that $\int_{\real^2}v_2 w\d x = 0$. 
We thus arrive at the following relation, 
\begin{equation}\label{Omegadef}
  \tilde \Omega \,=\, \tilde \Omega[w] \,:=\,
  -\frac{2}{d}\int_{\real^2}v_2(-x-2x_d)w(x)\d x~,
\end{equation}
which determines $\tilde \Omega$ as a function of $w$. 

As we shall see, since the vorticity profile $w$ is normalized 
by \eqref{wnorm}, the corresponding velocity field $v = K[w]$ 
satisfies
\[
  -v_2(-x-2x_d) \,\sim\, \frac{1}{2\pi d}~, \qquad \hbox{as} 
  \quad d \to \infty~,
\]
for any fixed $x \in \real^2$. In view of \eqref{Omegadef}, this means
that $\tilde \Omega \sim (\pi d^2)^{-1}$ as $d \to \infty$. 
In particular, if we take formally the limit $d \to \infty$ in 
\eqref{weq}, we obtain the limiting equation $v\cdot\nabla w = 0$, 
which simply means that $w$ (or $v$) is a stationary solution of 
Euler's equation. In what follows, we assume that the limiting
profile $w_*$ is {\em radially symmetric} and {\em stable} 
in the sense of Arnold \cite{Arn,MP94}. Roughly speaking, this 
means that $w_*(x)$ is a strictly decreasing function of $|x|$.
Given any such profile, we shall construct perturbatively a 
family of approximate solutions of \eqref{weq}, indexed by the
parameter $d \gg 1$, which converge to $w_*$ as $d \to \infty$.
Under natural assumptions, these approximate solutions are uniquely 
determined by the asymptotic profile $w_*$. 

An important question that we leave open here is whether we can
actually construct {\em exact} solutions of \eqref{weq} which converge
to $w_*$ as $d \to \infty$, in which case our approximate solutions
could be recovered by truncating the asymptotic expansion of the exact
solutions.  As was mentioned in the introduction, it should be
possible to prove the existence of such solutions, at least for a
particular class of compactly supported profiles $w_*$, by adapting
the variational techniques of D. Smets and J. van Schaftingen
\cite{SVS}. It might also be possible to construct exact solutions of
\eqref{weq} for more general profiles using a fixed point argument of
Nash-Moser type. We hope to clarify these issues in a future work.

\subsection{Asymptotic profile and functional setting}

Let $w_* \in \cS(\real^2)$ be a radially symmetric, nonnegative 
function satisfying the normalization condition \eqref{wnorm}, 
and let $v_* = K[w_*]$ be the velocity field obtained from
$w_*$ via the Biot-Savart law \eqref{BS}. We can thus write
\begin{equation}\label{qQdef}
  w_*(x) \,=\, \frac{1}{\pi}\,q(|x|^2)~, \qquad 
  v_*(x) \,=\, \frac{1}{2\pi} \frac{x^\perp}{|x|^2}\,Q(|x|^2)~,
  \qquad x \in \real^2~,
\end{equation}
where $q : [0,\infty) \to \real_+$ is a smooth, rapidly decreasing
function and $Q(r) = \int_0^r q(s)\d s$. Note that $Q(r) \to 1$ as 
$r \to \infty$. We assume that $w_*$ satisfies the following 
{\em strong stability conditions\:}
\begin{equation}\label{SS}
  q'(r) \,<\, 0 ~\hbox{for all}~ r \ge 0~, \qquad \hbox{and}\qquad
  \sup_{r > 0} \frac{-r^2 q'(r)}{Q(r)} < 1~.
\end{equation}
The first condition in \eqref{SS} implies of course that $q(r) > 0$
for all $r \ge 0$, so that $Q(r) > 0$ for all $r > 0$. In particular, 
compactly supported asymptotic profiles $w_*$ are excluded. This 
condition also implies that $w_*$ is a stable solution of the 
two-dimensional inviscid vorticity equation, with respect to 
perturbations in $L^1 \cap L^\infty$ \cite{MP85}. The second 
assumption in \eqref{SS} is more technical in nature, and can 
probably be relaxed. It is satisfied, for instance, if $q(r) = 
\gamma e^{-\gamma r}$ for some $\gamma > 0$. Note that we always 
have
\[
  \sup_{r > 0} \frac{-r^2 q'(r)}{Q(r)} > \frac14~.
\]
Indeed, if we assume on the contrary that $Q(r) + 4r^2 q'(r) \ge 0$ 
for all $r > 0$, then the function $h(r) = Q''(r) + Q(r)/(4r^2)$ 
is nonnegative and satisfies $h(r) \sim q(0)/(4r)$ as $r \to 0$, 
and $h(r) \sim 1/(4r^2)$ as $r \to \infty$. Thus $\sqrt{r}h \in 
L^1((0,\infty))$, but if we integrate by parts we obtain 
\[
  \int_0^\infty \sqrt{r} h(r)\d r \,=\, 
  \int_0^\infty \sqrt{r} \Bigl(Q''(r) + \frac{1}{4r^2}Q(r)\Bigr)\d r
  \,=\, 0~,
\] 
which yields a contradiction. Finally, in addition to \eqref{SS}, we 
also assume that $q^2/q'$ decays rapidly at infinity\:
\begin{equation}\label{qdecay}
  \sup_{r > 0} \frac{r^k q(r)^2}{|q'(r)|} \,<\, \infty~, \qquad
  \hbox{for all}~k \in \natural~.  
\end{equation}

As was already observed, the asymptotic profile $w_*$ is 
already an approximate solution of \eqref{weq}, \eqref{Omegadef}
in the sense that, if we substitute $(w_*,v_*)$ for $(w,v)$ 
in \eqref{weq}, the left-hand side converges to zero as $d \to 
\infty$. Our goal is to construct here more accurate approximations, 
which take into account the interaction of the vortices. We look
for solutions of the form
\begin{equation}\label{omudef}
  w \,=\, w_* + \omega~, \qquad v = v_* + u~,
\end{equation}
where $u = K[\omega]$ is the velocity field obtained from $\omega$ 
via the Biot-Savart law \eqref{BS}. The symmetry \eqref{wsym} implies
that
\begin{equation}\label{omsym}
  \omega(x_1,-x_2) \,=\, \omega(x_1,x_2)~, \quad
  u_1(x_1,-x_2) \,=\, -u_1(x_1,x_2)~, \quad
  u_2(x_1,-x_2) \,=\, u_2(x_1,x_2)~,
\end{equation}
and in agreement with \eqref{wnorm}, \eqref{wmom} we impose
\begin{equation}\label{omom}
  \int_{\real^2} \omega(x)\d x \,=\, \int_{\real^2} x_1 \omega(x)\d x 
  \,=\, 0~.
\end{equation}
Finally, we assume without loss of generality that $\omega$ has 
no radially symmetric component, namely
\begin{equation}\label{omzero}
  \int_0^{2\pi} \omega(r\cos\theta,r\sin\theta)
  \d \theta \,=\, 0 \qquad \hbox{for all} \quad r \ge 0~.
\end{equation}
We can always realize \eqref{omzero} by including, if necessary, the
radially symmetric part of $\omega$ into the asymptotic profile
$w_*$. In this respect, it is important to note that both conditions
in \eqref{SS} are open. 

Inserting \eqref{omudef} into \eqref{weq}, we obtain for $\omega$ the
following equation
\begin{equation}\label{omeq}
  \Lambda \omega + u\cdot\nabla\omega + R_d[\omega] \,=\, 0~,
\end{equation}
where $\Lambda$ is the linearized operator defined by
\begin{equation}\label{Lamdef}
  \Lambda \omega \,=\, v_* \cdot\nabla \omega + u\cdot
  \nabla w_*~,
\end{equation}
and $R_d[\omega]$ is a remainder term which depends on the 
distance $d$ between the vortex centers
\begin{equation}\label{Rdef}
  R_d[\omega](x) \,=\, \Bigl(v_*(x+2x_d)  - u(-x-2x_d) -\tilde 
  \Omega[w_*+\omega](x+x_d)^\perp\Bigr)\cdot\nabla 
  (w_*(x)+\omega(x))~.
\end{equation}
In \eqref{omeq}--\eqref{Rdef}, it is understood that $u = K[\omega]$ 
is the velocity field associated to $\omega$. 

We look for solutions $\omega$ of \eqref{omeq} in the Hilbert space
\begin{equation}\label{Xdef}
  X \,=\, \Bigl\{\omega \in L^2(\real^2)\,\Big|\, \int_{\real^2}
  |\omega(x)|^2 p(|x|^2)\d x < \infty\Bigr\}~,
\end{equation}
equipped with the scalar product
\[
  \langle \omega_1,\omega_2\rangle \,=\, \int_{\real^2}
  \omega_1(x)\omega_2(x)\,p(|x|^2)\d x~, \qquad \omega_1,\omega_2 
  \in X~,
\]
and with the associated norm $\|\omega\| \,=\,\langle \omega,\omega
\rangle^{1/2}$. Here the weight $p : [0,\infty) \to \real_+$ is defined by
\begin{equation}\label{pdef}
  p(r) \,=\, \frac{-1}{q'(r)}~, \qquad r \ge 0~.
\end{equation}
The reason for this particular choice is that the linear operator
$\Lambda$ has nice properties in the space $X$, see Section~\ref{Lamss} 
below. In view of \eqref{qdecay}, the asymptotic profile $w_*$ and
all its moments belongs to $X$\: for any $k \in \natural$, we have
\[
  \||x|^{2k}w_*\|^2 \,=\, \frac{1}{\pi^2}\int_{\real^2}|x|^{4k}
  q(|x|^2)^2 p(|x|^2)\d x \,=\, \frac{1}{\pi}\int_0^\infty
  r^{2k}q(r)^2p(r)\d r \,<\, \infty~.
\]

It is easy to verify that the operator $\Lambda$ commutes with the 
rotations about the origin in $\real^2$, see Lemma~\ref{Lampol} below.
It is thus natural to use polar coordinates $(r,\theta)$ in the plane, 
and to decompose our space $X$ as a direct sum
\begin{equation}\label{Xdec}
  X \,=\, \mathop{\oplus}\limits_{n=0}^\infty X_n \,=\, 
  \mathop{\oplus}\limits_{n=0}^\infty P_n X~,
\end{equation}
where $P_n$ is the orthogonal projection in $X$ defined by the formula
\[
  (P_n \omega)(r\cos\theta,r\sin\theta) \,=\, \frac{2-\delta_{n,0}}{2\pi}
  \int_0^{2\pi} \omega(r\cos\theta',r\sin\theta')\cos(n(\theta-\theta'))
  \d\theta'~, \quad n \in \natural~.
\]
In particular, $X_0 = P_0 X$ is the subspace of all radially symmetric 
functions, and for $n \ge 1$ the subspace $X_n = P_n X$ contains 
functions of the form $\omega(r\cos\theta,r\sin\theta) = 
a_1(r)\cos(n\theta) + a_2(r)\sin(n\theta)$. With this notation, 
condition \eqref{omzero} means that $P_0 \omega = 0$. 

\subsection{The linearized operator and its right-inverse}\label{Lamss}

We now discuss the main properties of the linearized operator
$\Lambda$ defined in \eqref{Lamdef}. In the particular case where
$w_*$ is the profile $G$ of Oseen's vortex \eqref{Gdef}, the operator
$\Lambda$ was studied in detail in \cite{GW05,Ma1}, and we shall
obtain here analogous results in a more general situation.

From \eqref{Lamdef} we know that $\Lambda = \Lambda_1 + \Lambda_2$, 
where $\Lambda_1 \omega = v_*\cdot\nabla\omega$ and $\Lambda_2 
\omega = K[\omega]\cdot\nabla w_*$. As is easily verified, 
$\Lambda_2$ is compact in $X$, while $\Lambda_1$ is unbounded. 
The maximal domain of $\Lambda$ is therefore
\begin{equation}\label{Lamdom}
  D(\Lambda) \,=\, D(\Lambda_1) \,=\, \{\omega \in X \,|\,
  v_*\cdot\nabla \omega \in X\}~.
\end{equation}
The most remarkable property of this operator is that it is 
{\em skew-symmetric} in $X$. 

\begin{lemma}\label{skew}
For all $\omega_1, \omega_2 \in D(\Lambda)$, we have $\langle 
\Lambda \omega_1,\omega_2\rangle + \langle \omega_1,\Lambda
\omega_2\rangle = 0$. 
\end{lemma}

\proof We shall prove in fact that both operators $\Lambda_1, 
\Lambda_2$ are skew-symmetric. First, since the weight $p(|x|^2)$ is 
radially symmetric, we have
\[
  \langle v_*\cdot\nabla\omega_1,\omega_2\rangle + \langle 
  \omega_1,v_* \cdot \nabla\omega_2\rangle \,=\, 
  \int_{\real^2} p(|x|^2)\,v_*\cdot\nabla(\omega_1\omega_2)\d x
  \,=\, 0~, 
\]
because the velocity field $p(|x|^2)v_*(x)$ is divergence-free. 
Next, since $\nabla w_*(x) = (2/\pi)x q'(|x|^2)$, we have
\[
  \langle u_1\cdot\nabla w_*,\omega_2\rangle + \langle 
  \omega_1,u_2 \cdot \nabla w_*\rangle \,=\, -\frac{2}{\pi}
  \int_{\real^2}\Bigl((x\cdot u_1)\omega_2 + (x\cdot u_2)\omega_1
  \Bigr)\d x \,=\, 0~, 
\]
see e.g. \cite[Lemma~4.8]{GW05}. Combining both equalities, 
we obtain the desired result. \QED

\medskip
As in the Gaussian case \cite{GW05}, the operator $\Lambda$ 
is invariant under rotations about the origin in the plane 
$\real^2$. It is thus natural to work in polar coordinates 
$(r,\theta)$, and to develop the vorticity $\omega(r\cos\theta,
r\sin\theta)$ in Fourier series with respect to the angular 
variable $\theta$. In these variables, the action of $\Lambda$ 
can be described fairly explicitly. Let
\begin{equation}\label{phidef}
  \phi(r) \,=\, \frac{Q(r^2)}{2\pi r^2}~, \qquad
  g(r) \,=\, -\frac{2 q'(r^2)}{\pi}~, \qquad r > 0~.
\end{equation}
Then we have the following result\:

\begin{lemma}\label{Lampol}
Fix $n \in \natural$. If $\omega = a_n(r)\sin(n\theta)$, then 
$\Lambda \omega = n[\phi(r)a_n(r) - g(r)A_n(r)]\cos(n\theta)$, 
where $A_n$ is the regular solution of the differential 
equation
\begin{equation}\label{Andef}
  -A_n''(r) - \frac1r A_n'(r) + \frac{n^2}{r^2}A_n(r) \,=\, 
  a_n(r)~, \qquad r > 0~.
\end{equation}
Similarly, if $\omega = -a_n(r)\cos(n\theta)$, then 
$\Lambda \omega = n[\phi(r)a_n(r) - g(r)A_n(r)]\sin(n\theta)$. 
\end{lemma}

\proof
If $n = 0$, namely if $\omega$ is radially symmetric, it is 
straightforward to verify that $\Lambda \omega = 0$. Thus 
we assume that $n \ge 1$, and that $\omega = a_n(r)\sin(n\theta)$. 
If $\psi$ denotes the stream function defined by $-\Delta \psi = 
\omega$, we have $\psi = A_n(r)\sin(n\theta)$, where $A_n$ is the 
regular solution of \eqref{Andef}, namely
\[
  A_n(r) \,=\, \frac{1}{2n}\left(\int_0^r \Bigl(\frac{s}{r}\Bigr)^n
  sa_n(s)\d s + \int_r^\infty \Bigl(\frac{r}{s}\Bigr)^n
  sa_n(s)\d s\right)~, \qquad r > 0~. 
\]
The velocity field $u = K[\omega] = -\nabla^\perp\psi$ is thus
\[
  u \,=\, \frac{n}{r}A_n(r)\cos(n\theta)\mathbf{e}_r - A_n'(r)
  \sin(n\theta)\mathbf{e}_\theta~,   \qquad \hbox{where}\quad
  \mathbf{e}_r \,=\, \frac{x}{|x|}~, \quad \mathbf{e}_\theta \,=\, 
  \frac{x^\perp}{|x|}~.
\]
Since $\Lambda \omega = v_*\cdot\nabla \omega + u\cdot\nabla w_*$, 
we conclude that
\[
  \Lambda \omega \,=\, \frac{Q(r^2)}{2\pi r^2}\,n a_n(r)
  \cos(n\theta) + \frac{n}{r}A_n(r)\cos(n\theta)\,\frac{2r}{\pi}
  q'(r^2)~,
\]
which, in view of \eqref{phidef}, is the desired result.  
The case where $\omega = -a_n(r)\cos(n\theta)$ is similar. 
\QED

\medskip
As an application of Lemma~\ref{Lampol}, we can characterize 
the kernel of the operator $\Lambda$. We already know that
$\Lambda\omega = 0$ if $\omega$ is radially symmetric. Moreover, 
differentiating the identity $v_* \cdot\nabla w_* = 0$ with 
respect to $x_1$ and $x_2$, we obtain $\Lambda(\partial_1 w_*)
= \Lambda(\partial_2 w_*) = 0$. As in the Gaussian case 
\cite{Ma1}, we conclude\:

\begin{lemma}\label{Lamker}
$\Ker(\Lambda) \,=\, X_0 \oplus\{\alpha_1 \partial_1w_* + 
\alpha_2 \partial_2 w_*\,|\, \alpha_1,\alpha_2 \in \real\}$. 
\end{lemma}

\proof
Since the decomposition \eqref{Xdec} is invariant under the action 
of $\Lambda$, it is sufficient to characterize the kernel in each 
subspace $X_n$. The case $n = 0$ is trivial, because $X_0 \subset
\Ker(\Lambda)$, hence we assume from now on that $n \ge 1$. If
$\omega = a_n(r)\sin(n\theta)$ satisfies $\Lambda \omega = 0$, 
we know from Lemma~\ref{Lampol} that $\phi a_n - g A_n = 0$. 
In view of \eqref{Andef}, this can be written in the equivalent 
form
\begin{equation}\label{Aneq}
  -A_n''(r) - \frac1r A_n'(r) + \Bigl(\frac{n^2}{r^2} - 
  \frac{g(r)}{\phi(r)}\Bigr)A_n(r) \,=\, 0~, \qquad r > 0~.
\end{equation}
Now, the second assumption in \eqref{SS} means that
\[
  \sup_{r > 0}\frac{r^2g(r)}{\phi(r)} \,=\, 
  \sup_{r > 0}\frac{4 r^2 |q'(r)|}{Q(r)} \,<\, 4~. 
\]
Thus, if $n \ge 2$, the ``potential'' term $(n^2/r^2 - g/\phi)$ 
in \eqref{Aneq} is positive, and since $A_n(r) \to 0$ as $r \to 0$ 
and $r \to \infty$, the maximum principle implies that $A_n = 0$, 
hence also $a_n = 0$. Thus $\Ker(\Lambda) \cap X_n = \{0\}$ 
if $n \ge 2$. 

In the particular case $n = 1$, it is easy to verify that 
$A_1(r) = r\phi(r)$ is the regular solution of \eqref{Aneq}. 
Using \eqref{Andef} we find $a_1(r) = rg(r)$, so that $\omega = 
a_1(r)\sin(\theta) = -\partial_2 w_*$. Similarly, $a_1(r)\cos(\theta) 
= -\partial_1 w_*$, hence the kernel of $\Lambda$ in $X_1$ is 
spanned by the functions $\{\partial_1 w_*,\partial_2 w_*\}$. 
\QED

\medskip
Using the same arguments as in \cite{Ma1}, one can show that the 
operator $\Lambda$ is not only skew-symmetric, but also 
{\em skew-adjoint} in $X$. This implies that $\Ker(\Lambda)
= \Ran(\Lambda)^\perp$, hence
\[
  \overline{\Ran(\Lambda)} \,=\, \Ker(\Lambda)^\perp~.
\]
Let 
\begin{equation}\label{Ydef}
  Y \,=\, \{\omega \in X\,|\, |x|^2 \omega \in X\}~.
\end{equation}
We now show that $\Ker(\Lambda)^\perp \cap Y \subset \Ran(\Lambda)$, 
and we establish a semi-explicit formula for the inverse of 
$\Lambda$ on that subspace. 

\begin{proposition}\label{Laminv}
If $f \in X_n \cap Y$ for some $n \ge 2$, there exists a unique 
$\omega \in X_n \cap D(\Lambda)$ such that $\Lambda\omega = f$. 
Specifically, if $f = b_n(r)\cos(n\theta)$, then $\omega = a_n(r)
\sin(n\theta)$, where
\begin{equation}\label{anAn}
  a_n(r) \,=\, \frac{g(r)}{\phi(r)}A_n(r) + \frac{b_n(r)}{n\phi(r)}~,
\end{equation}
and $A_n$ is the regular solution of the differential equation
\begin{equation}\label{Anbn}
  -A_n''(r) - \frac1r A_n'(r) + \Bigl(\frac{n^2}{r^2} - 
  \frac{g(r)}{\phi(r)}\Bigr)A_n(r) \,=\, \frac{b_n(r)}{n\phi(r)}~, 
  \qquad r > 0~.
\end{equation}
Similarly, if $f = b_n(r)\sin(n\theta)$, then $\omega = -a_n(r)
\cos(n\theta)$. 
\end{proposition}

\proof If $\omega = a_n(r)\sin(n\theta)$, then $\Lambda\omega =  
n[\phi(r)a_n(r) - g(r)A_n(r)]\cos(n\theta)$ by Lemma~\ref{Lampol}, 
where $A_n$ satisfies \eqref{Andef}. The equation we have to solve 
is therefore $n(\phi a_n - g A_n) = b_n$, which gives \eqref{anAn}. 
Moreover, combining \eqref{anAn} and \eqref{Andef}, we obtain 
\eqref{Anbn}. 

Proceeding as in \cite[Lemma~3.4]{Ga}, we now show that \eqref{Anbn}
has a unique regular solution, and we establish a representation
formula. As we observed in the proof of Lemma~\ref{Lamker}, the
``potential'' term $(n^2/r^2 - g/\phi)$ in \eqref{Aneq} is positive if
$n \ge 2$.  Let $\psi_+$, $\psi_-$ be the (unique) solutions of the
homogeneous equation \eqref{Aneq} such that
\begin{equation}\label{psipmdef}
  \psi_-(r) \,\sim\, r^n \quad \hbox{as }r \to 0~, \quad 
  \hbox{and}\quad \psi_+(r) \,\sim\, r^{-n} \quad \hbox{as }r 
  \to \infty~.
\end{equation}
By the maximum principle, the functions $\psi_+$, $\psi_-$ are strictly 
monotone and linearly independent. The Wronskian determinant 
$W = \psi_+ \psi_-' - \psi_- \psi_+'$ satisfies $W' + W/r = 0$, 
hence $W(r) = 2n\kappa/r$ for some $\kappa > 0$, and we also have
\[
  \psi_-(r) \,\sim\, \kappa r^n \quad \hbox{as }r \to \infty~, \quad 
  \hbox{and}\quad \psi_+(r) \,\sim\, \kappa r^{-n} \quad \hbox{as }r 
  \to 0~.
\]
With these notations, the unique regular solution of \eqref{Anbn}
has the following expression\:
\begin{equation}\label{Anrep}
  A_n(r) \,=\, \psi_+(r)\int_0^r \frac{\psi_-(s)}{W(s)}\,\frac{b_n(s)}
  {n\phi(s)}\d s + \psi_-(r)\int_r^\infty \frac{\psi_+(s)}{W(s)}
  \,\frac{b_n(s)}{n\phi(s)}\d s~, \quad r > 0~. 
\end{equation}

If $f = b_n(r)\cos(n\theta) \in X_n$, it is straightforward
to verify that the function $A_n$ defined by \eqref{Anrep}
is continuous and vanishes at the origin and at infinity. 
Moreover, we know from \eqref{phidef} that $\phi(r) \sim 
1/(2\pi r^2)$ as $r \to \infty$. Thus, if we assume that 
$f \in X_n \cap Y$, we see that the function $a_n$ defined
by \eqref{anAn} satisfies $\int_0^\infty a_n(r)^2p(r^2)r\d r < \infty$.
As a consequence, if $\omega = a_n(r)\sin(n\theta)$, we conclude 
that $\omega \in X_n \cap D(\Lambda)$, and $\Lambda\omega = f$ 
by construction.
\QED

\begin{remark}\label{none}
If $n = 1$, the conclusion of Proposition~\ref{Laminv} fails because 
$\partial_j w_* \in X_1 \cap \Ker(\Lambda)$ for $j = 1,2$. However, 
if $f \in X_1 \cap Y$ satisfies $\langle f,\partial_jw_*\rangle = 0$ 
for $j = 1,2$, one can show that there exists a unique $\omega \in X_1 
\cap D(\Lambda) \cap \Ker(\Lambda)^\perp$ such that $\Lambda \omega = f$.
\end{remark}

\subsection{The perturbation expansion}

Equipped with the technical results of the previous section, 
we now go back to equation \eqref{omeq}, which we want to 
solve perturbatively for large $d$. This equation can be written
as $\Lambda\omega + N_d[\omega] = 0$, where $N_d[\omega] = 
u\cdot\nabla\omega + R_d[\omega]$. Before starting the calculations, 
we briefly explain why we expect to find a unique solution, 
under our symmetry assumptions. 

First, if $\omega \in X$ satisfies \eqref{omsym}--\eqref{omzero}, 
then $\omega \in \Ker(\Lambda)^\perp$, hence $\omega$ is uniquely
determined by $\Lambda\omega$. Indeed, as was already observed, 
\eqref{omzero} means that $P_0 \omega = 0$. Moreover, it follows 
from \eqref{omsym}, \eqref{omom} that
\[
  \langle \partial_j w_*,\omega\rangle \,=\, -\frac{2}{\pi}
  \int_{\real^2} x_j \omega \d x \,=\, 0~, \qquad j = 1,2~,
\]
hence $\omega \in \Ker(\Lambda)^\perp$ by Lemma~\ref{Lamker}.  Next,
if $\omega$ and $u = K[\omega]$ have the symmetries \eqref{omsym}, it
is straightforward to verify that the nonlinearity in \eqref{omeq}
satisfies $N_d[\omega](x_1,-x_2) = -N_d[\omega](x_1,x_2)$ for all $x =
(x_1,x_2) \in \real^2$. This implies that $P_0 N_d[\omega] = 0$ and
$\langle \partial_1 w_*,N_d[\omega]\rangle = 0$. Moreover, we have
$\langle \partial_2 w_*,N_d[\omega]\rangle = 0$ by construction,
because this is the relation we imposed to determine the angular speed
$\tilde \Omega$ in \eqref{Omegadef}. Thus, we see that $N_d[\omega]
\in \Ker(\Lambda)^\perp$, and if we can prove in addition that $|x|^2
N_d[\omega] \in X$, then Proposition~\ref{Laminv} (and
Remark~\ref{none}) will imply that $N_d[\omega] \in \Ran(\Lambda)$.
We can therefore hope to find a unique $\omega \in \Ker(\Lambda)^\perp
\cap D(\Lambda)$ such that $\Lambda\omega + N_d[\omega] = 0$.

To begin our perturbative approach, we compute the remainder 
term \eqref{Rdef} for $\omega = 0$, namely
\begin{equation}\label{Rddef}
  R_d(x) \,\equiv\, R_d[0](x) \,=\, \Bigl(v_*(x+2x_d) - \tilde \Omega[w_*]
  (x+x_d)^\perp \Bigr)\cdot\nabla w_*(x)~, \qquad x \in \real^2~.
\end{equation}
From \eqref{qQdef}, we know that
\[
  v_*(x) \,=\, \frac{1}{2\pi}\frac{x^\perp}{|x|^2}(1 - \tilde Q(|x|^2))~,
  \qquad \hbox{where}\quad \tilde Q(r) = \int_r^\infty q(s)\d s~.
\]
By assumption, the term $\tilde Q(|x|^2)$ decays faster than any 
inverse power of $|x|$ as $|x| \to \infty$, hence we can neglect its
contribution in our calculations. For any fixed $x \in \real^2$ we 
thus have
\begin{equation}\label{vstardec}
  v_*(x+2x_d) \,=\, \frac{1}{2\pi}\frac{(2x_d)^\perp}{|2x_d|^2}
  + \frac{1}{2\pi} V(x,2x_d) + \cO\Bigl(\frac{1}{d^\infty}\Bigr)~,
  \qquad \hbox{as}\quad d \to \infty~,
\end{equation}
where
\[
  V(x,y) \,=\, \frac{(x+y)^\perp}{|x+y|^2} - \frac{y^\perp}{|y|^2}~.
\]
Setting $x = (r\cos\theta,r\sin\theta)$, $y = 2x_d = (d,0)$, and 
proceeding as in \cite[Lemma~3.2]{Ga}, we find
\begin{equation}\label{Vexp}
  V_1(x,y) \,=\, \frac{1}{d}\sum_{n=1}^\infty (-1)^n\,\frac{r^n}{d^n}
  \,\sin(n\theta)~, \qquad
  V_2(x,y) \,=\, \frac{1}{d}\sum_{n=1}^\infty (-1)^n\,\frac{r^n}{d^n}
  \,\cos(n\theta)~.
\end{equation}
In particular, returning to \eqref{vstardec} and using definition 
\eqref{Omegadef}, we obtain
\begin{equation}\label{Omexp}
  \tilde \Omega[w_*] \,=\, \frac{2}{d}\int_{\real^2}(v_*)_2(x+2x_d)w_*(x)
  \d x \,=\, \frac{1}{\pi d^2} + \cO\Bigl(\frac{1}{d^\infty}\Bigr)~,
  \qquad \hbox{as}\quad d \to \infty~.
\end{equation}
Note that the term $V(x,2x_d)$ in \eqref{vstardec} gives no 
contribution to the angular velocity $\tilde \Omega[w_*]$. On 
the other hand, inserting \eqref{vstardec}, \eqref{Omexp} into 
\eqref{Rddef} and using the expansion \eqref{Vexp} together with
the relation $\nabla w_* = -xg(|x|)$, where $g$ is defined 
in \eqref{phidef}, we find for $x = (r\cos\theta,r\sin\theta)$\:
\begin{equation}\label{Rdexp}
  R_d(x) \,=\, \frac{g(r)}{2\pi} \sum_{n=2}^\infty (-1)^n\,
  \frac{r^n}{d^n}\,\sin(n\theta) + \cO\Bigl(\frac{1}{d^\infty}\Bigr)~,
  \qquad \hbox{as}\quad d \to \infty~.
\end{equation}

Motivated by this result, we now construct inductively an 
approximate solution of \eqref{omeq} of the form
\begin{equation}\label{omapp}
  \omega(x) \,=\, \sum_{n=2}^\ell \frac{1}{d^n}\,\omega^{(n)}(x)~, \qquad
  u(x) \,=\, \sum_{n=2}^\ell \frac{1}{d^n}\,u^{(n)}(x)~,
\end{equation}
where each velocity profile $u^{(n)}$ is obtained from $\omega^{(n)}$
via the Biot-Savart law \eqref{BS}. The order $\ell$ of the
approximation is in principle arbitrary, but the complexity 
of the calculations increases rapidly with $\ell$, and we 
shall restrict ourselves to $\ell = 4$ for simplicity. 
Of course, we assume that the symmetry and normalization conditions
\eqref{omsym}--\eqref{omzero} hold at each order of the approximation. 
In particular, we have
\begin{equation}\label{omjmom}
  \int_{\real^2} \omega^{(n)}(x)\d x \,=\, 
  \int_{\real^2} x_1\omega^{(n)}(x)\d x \,=\, 
  \int_{\real^2} x_2\omega^{(n)}(x)\d x \,=\, 0~, 
\end{equation}
for all $n \in \{2,\dots,\ell\}$. In view of \cite[Appendix~B]{GW02}, 
this implies that the velocity field $u^{(n)}(x)$ decays at least 
as fast as $|x|^{-3}$ when $|x| \to \infty$. It follows that the 
term $u(-x-2x_d)$ in \eqref{Rdef} is $\cO(d^{-5})$ as $d \to \infty$, 
and will therefore not contribute to $\omega^{(n)}$ for $n \le 4$. 
For the same reason,
\[
  \tilde \Omega[w_* + \omega] \,=\, \frac{2}{d}\int_{\real^2}
  \Bigl((v_*)_2(x+2x_d) - u_2(-x-2x_d)\Bigr)\Bigl(w_*(x) + 
  \omega(x)\Bigr)\d x \,=\, \frac{1}{\pi d^2} + \cO\Bigl(\frac{1}
  {d^6}\Bigr)~,
\]
as $d \to \infty$. Indeed, the leading term $\tilde\Omega[w_*]$ 
was computed in \eqref{Omexp}, and we know that the contribution 
of $u_2(-x-2x_d)$ is negligible. Moreover, using \eqref{vstardec}, 
\eqref{Vexp}, and \eqref{omjmom}, it is easy to verify that 
$\int (v_*)_2(x+2x_d) \omega(x)\d x = \cO(d^{-5})$, as $d \to 
\infty$. Summarizing, we have shown that
\begin{align}\nonumber
  R_d[\omega](x) \,&=\, \Bigl(v_*(x+2x_d)  -\tilde \Omega[w_*]
  (x+x_d)^\perp\Bigr)\cdot\nabla (w_*(x)+\omega(x)) + \cO\Bigl(\frac{1}
  {d^5}\Bigr) \\ \label{Rdtrunc}
  \,&=\, R_d(x) + \Bigl(\frac{1}{2\pi}V(x,2x_d) - \frac{x^\perp}{\pi
  d^2}\Bigr)\cdot\nabla \omega(x) + \cO\Bigl(\frac{1}{d^5}\Bigr)~,  
\end{align}
as $d \to \infty$. Similarly, the quadratic term $u\cdot\nabla\omega$
in \eqref{omeq} satisfies
\begin{equation}\label{quadtrunc}
  u\cdot\nabla\omega \,=\, \frac{1}{d^4} u^{(2)}\cdot\nabla
  \omega^{(2)} + \cO\Bigl(\frac{1}{d^5}\Bigr)~, \qquad \hbox{as}
  \quad d \to \infty~.
\end{equation}
 
It is now a straightforward task to determine the first vorticity 
profiles in the expansion \eqref{omapp}. From \eqref{Rdtrunc}, 
\eqref{quadtrunc}, we know that the nonlinearity $N_d[\omega] = 
u\cdot\nabla\omega + R_d[\omega]$ in \eqref{omeq} satisfies, 
for $x = (r\cos\theta,r\sin\theta) \in \real^2$, 
\begin{equation}\label{Ndexp1}
  N_d[\omega](x) \,=\, \frac{g(r)}{2\pi}\left(\frac{r^2}{d^2}
  \sin(2\theta) - \frac{r^3}{d^3}\sin(3\theta)\right) + 
  \cO\Bigl(\frac{1}{d^4}\Bigr)~, \qquad \hbox{as}
  \quad d \to \infty~.
\end{equation}
Thus, to ensure that $\Lambda\omega + N_d[\omega] = \cO(d^{-4})$, we 
must impose
\begin{equation}\label{om23}
  \Lambda \omega^{(n)} + \frac{g(r)}{2\pi}(-1)^n r^n \sin(n\theta)
  \,=\, 0~, \qquad \hbox{for} \quad n = 2,3~.
\end{equation}
By Proposition~\ref{Laminv}, Eq.~\eqref{om23} has a unique solution
$\omega^{(n)} \in X_n \cap D(\Lambda)$ of the form
\begin{equation}\label{omjexp}
  \omega^{(n)}(x) \,=\, a_n(r)\cos(n\theta)~, \qquad
  u^{(n)}(x) \,=\, -\frac{n}{r}A_n(r)\sin(n\theta)\mathbf{e}_r 
  - A_n'(r)\cos(n\theta)\mathbf{e}_\theta~,
\end{equation}
where $a_n(r), A_n(r)$ are given by \eqref{anAn}, \eqref{Anbn} with
$b_n(r) = (-1)^nr^ng(r)/(2\pi)$. As is easily verified, the symmetry
conditions \eqref{omsym}--\eqref{omzero}, are satisfied by the
profiles $\omega^{(n)}, u^{(n)}$ for $n = 2,3$, and the velocity
$|u^{(n)}(x)|$ decays like $|x|^{-n-1}$ as $|x| \to \infty$.

Computing the profiles $\omega^{(4)}, u^{(4)}$ is more cumbersome, but 
also more representative of what happens in the general case. First
of all, the quadratic term \eqref{quadtrunc} is no longer negligible, 
and using \eqref{omjexp} for $n = 2$ we find
\[
  u^{(2)}(x)\cdot\nabla\omega^{(2)}(x) \,=\, B_1(r)\sin(4\theta)~, 
  \qquad \hbox{where}\quad B_1(r) \,=\, \frac1r \Bigl(
  A_2'(r)a_2(r) - A_2(r)a_2'(r)\Bigr)~.   
\]
Note that, to ensure that $u^{(2)}\cdot\nabla\omega^{(2)} \in X$, 
we need an assumption on the second derivative of the function 
$q$ appearing in \eqref{qQdef}. For instance, in analogy with 
\eqref{qdecay}, one can impose
\begin{equation}\label{q2decay}
  \sup_{r > 0} \frac{r^k q''(r)^2}{|q'(r)|} \,<\, \infty~, \qquad
  \hbox{for all}~k \in \natural~.  
\end{equation}
Next, we must compute the contribution of $\omega^{(2)}$ to 
the right-hand side of \eqref{Rdtrunc}. Using \eqref{Vexp}, 
\eqref{omjexp}, we obtain
\[
  \Bigl(\frac{d^2}{2\pi}V(x,2x_d) - \frac{x^\perp}{\pi
  }\Bigr)\cdot\nabla \omega^{(2)}(x) \,=\, B_2(r)\sin(4\theta) + 
  C(r)\sin(2\theta)~,
\]
where
\[
  B_2(r) \,=\, \frac{1}{4\pi}(2a_2(r)-ra_2'(r))~, \qquad
  C(r) \,=\, \frac{2}{\pi}a_2(r)~.
\]
Combining these results, we find instead of \eqref{Ndexp1}\:
\[
  N_d[\omega](x) \,=\, \frac{g(r)}{2\pi}\sum_{n=2}^3 (-1)^n
  \frac{r^n}{d^n}\sin(n\theta) + \frac{B(r)}{d^4}
  \sin(4\theta) + \frac{C(r)}{d^4}\sin(2\theta) + 
  \cO\Bigl(\frac{1}{d^5}\Bigr)~, 
\]
as $d \to \infty$, where $B(r) = B_1(r) + B_2(r) + r^4g(r)/(2\pi)$. 
Therefore, in addition to \eqref{om23}, we must impose 
\begin{equation}\label{om4}
  \Lambda \omega^{(4)} + B(r)\sin(4\theta) + C(r)\sin(2\theta) \,=\, 
  0~. 
\end{equation}
Using again Proposition~\ref{Laminv}, we see that \eqref{om4} 
has a unique solution $\omega^{(4)} \in (X_4 + X_2) \cap D(\Lambda)$
of the form $\omega^{(4)}(x) \,=\, a_4(r)\cos(4\theta) + \tilde 
a_2(r)\cos(2\theta)$, where $a_4(r)$ is given by \eqref{anAn}, 
\eqref{Anbn} with $n = 4$ and $b_4(r) = B(r)$, while $\tilde a_2(r)$ 
is given by the same relations with $n = 2$ and $b_2(r) = C(r)$. 
An explicit expression of the velocity profile $u^{(4)}$ can 
also be obtained, as in \eqref{omjexp}. 

\smallskip
Summarizing, we have shown\:

\begin{proposition}\label{pexp}
Let $w_*$ be a radially symmetric vorticity profile of the 
form \eqref{qQdef}, where the function $q$ satisfies 
\eqref{SS}, \eqref{qdecay}, \eqref{q2decay}, and let
\begin{equation}\label{wvasym}
  w(x) \,=\, w_*(x) + \sum_{n=2}^4\frac{1}{d^n}\omega^{(n)}(x)~, 
  \qquad
  v(x) \,=\, v_*(x) + \sum_{n=2}^4\frac{1}{d^n}u^{(n)}(x)~, 
\end{equation}
where the vorticity profiles $\omega^{(n)} \in X \cap D(\Lambda)
\cap \Ker(\Lambda)^\perp$ satisfy \eqref{om23}, \eqref{om4}, and 
the velocity profiles $u^{(n)}$ are obtained by the Biot-Savart law 
\eqref{BS}. Then $w$ is an asymptotic solution of Eqs.~\eqref{weq}, 
\eqref{Omegadef} in the sense that
\begin{equation}\label{wvapprox}
  \Bigl(v(x) - v(-x-2x_d) - \tilde\Omega[w](x+x_d)^\perp\Bigr)
  \cdot \nabla w(x) \,=\, \cO\Bigl(\frac{1}{d^5}\Bigr)~,
\end{equation}
in the topology of $X$ and uniformly on $\real^2$, as $d \to \infty$. 
\end{proposition}

The asymptotic expansion \eqref{omapp} is very natural, and it is
clear that it can be performed to any finite order $\ell \in \natural$
if we make appropriate assumptions on the derivatives of the profile
$q$, as in \eqref{q2decay}.  As was already mentioned, we also
conjecture that there exists an {\em exact} solution of \eqref{weq} for
$d \gg 1$ which coincides with \eqref{wvasym} up to corrections of
order $\cO(d^{-5})$.  It is important to notice that, under the
symmetry and normalization conditions \eqref{wnorm}--\eqref{wmom}, the
exact solution (if it exists) and the asymptotic expansion
\eqref{wvasym} are uniquely determined by the limiting profile $w_*$.

\begin{remark}\label{moregen}
Rather strong assumptions on the limiting profile $w_*$ were made in
this section to ensure that the asymptotic expansion \eqref{wvasym}
holds in the function space $X$ defined in \eqref{Xdef}, which we
believe is naturally associated to the problem. This does not restrict
the scope of our results here, because these conditions are
automatically fulfilled by the Gaussian profiles created by the
Navier-Stokes evolution.  However, within the framework of Euler's
equation, it is certainly interesting to construct interacting vortex
pairs with more general profiles, including compactly supported ones.
If we do not insist on controlling our expansion in the space $X$, the
calculations presented in this section show that the assumptions on
the function $q$ can be considerably relaxed.  The most important
point is that Eq.~\eqref{Anbn} should have a unique solution for $n
\ge 2$, and for $n = 1$ if the right-hand side satisfies some
orthogonality conditions.  This is definitely the case if the second
inequality in \eqref{SS} holds, but that condition does not imply that 
$q$ is strictly decreasing and can well be satisfied if $q$ is 
compactly supported.
\end{remark}

To conclude this section, we indicate how approximate solutions
of \eqref{Vrot} of the form \eqref{Vpair} can be obtained from 
Proposition~\ref{pexp} by a simple rescaling. Given $\alpha,
d > 0$, we consider the situation described in \eqref{twov}
with $\alpha_1 = \alpha_2 = \alpha$ and $r_1 = r_2 = d/2$.
If $w_*$ is a radially symmetric vorticity profile satisfying
the assumptions of Proposition~\ref{pexp}, we define, for 
all sufficiently small $\epsilon > 0$,
\begin{equation}\label{wvasym2}
  w_\epsilon(x) \,=\, w_*(x) + \sum_{n=2}^4\frac{\epsilon^n}{d^n}
  \,\omega^{(n)}(x)~, \qquad
  v_\epsilon(x) \,=\, v_*(x) + \sum_{n=2}^4\frac{\epsilon^n}{d^n}
  \,u^{(n)}(x)~, 
\end{equation}
where $\omega^{(n)}, u^{(n)}$ are as in \eqref{wvasym}. Then, 
by construction, the vorticity distribution
\begin{equation}\label{Vpair2}
  \omega_\epsilon(x) \,=\, \frac{\alpha}{\epsilon^2} \,w_\epsilon
   \Bigl(\frac{x-x_d}{\epsilon}\Bigr) + \frac{\alpha}{\epsilon^2} 
   \,w_\epsilon \Bigl(\frac{-x-x_d}{\epsilon}\Bigr)~,
\end{equation}
where $x_d = (d/2,0)$, is an approximate solution of Eq.~\eqref{Vrot} 
with $\Omega = \alpha/(\pi d^2)$. More precisely, it follows from 
\eqref{wvapprox} that
\begin{equation}\label{omuapprox}
  \partial_t \omega_\epsilon + (u_\epsilon -\Omega x^\perp) \cdot \nabla 
  \omega_\epsilon \,=\,  \cO(\epsilon)~, \qquad \hbox{as} \quad \epsilon 
  \to 0~,
\end{equation}
where $u_\epsilon$ is obtained from $\omega_\epsilon$ via the 
Biot-Savart law \eqref{BS}.

\section{Inviscid approximation of viscous vortex pairs}
\label{sec3}

In this final section, we describe in some detail the result of
\cite{Ga} in the particular case of a vortex pair, and we interprete
it using the approximate solutions of Euler's equation constructed
in Section~\ref{sec2}. Given $\alpha,d > 0$, we set $\Omega =
\alpha/(\pi d^2)$ and we denote by $\omega_0$ the vorticity
distribution \eqref{twov} where $\alpha_1 = \alpha_2 = \alpha$ and
$r_1 = r_2 = d/2$. For any $\nu > 0$, we consider the (unique) solution
$\omega^\nu(x,t)$ of the rotating viscous vorticity equation
\begin{equation}\label{Vrotvis}
  \partial_t \omega + (u - \Omega x^\perp) \cdot \nabla 
  \omega \,=\, \nu\Delta\omega~, \qquad x \in \real^2~, \quad
  t > 0~,
\end{equation}
with initial data $\omega_0$. Up to a rotation of angle $\Omega t$, 
the vorticity distribution $\omega^\nu(x,t)$ coincides with the solution
of the nonrotating equation \eqref{Veq} with the same initial data, 
which is studied in \cite{Ga}. The advantage of using a rotating 
frame is that the vortex centers remain fixed, instead of evolving
according to the point vortex dynamics \eqref{PW}. As a matter of
fact, Theorem~2.1 in \cite{Ga} establishes that $\omega^\nu(\cdot,t) 
\weakto \omega_0$ as $\nu \to 0$, for any $t > 0$. 

To obtain a more precise convergence result, we decompose 
the solution of \eqref{Vrotvis} into a sum of viscous 
vortices\:
\begin{equation}\label{omdecomp}
  \omega^\nu(x,t) \,=\, \frac{\alpha}{\nu t}\,w_1^\nu\Bigl(
  \frac{x-x_1}{\sqrt{\nu t}},t\Bigr) + \frac{\alpha}{\nu t}
  \,w_2^\nu\Bigl(\frac{x-x_2}{\sqrt{\nu t}},t\Bigr)~,
\end{equation}
where $x_1 = -x_2 = (d/2,0)$. As is shown in \cite{Ga}, both vorticity
profiles $w_1^\nu(\xi,t),w_2^\nu(\xi,t)$ can be approximated by the
same function
\begin{equation}\label{wappdef}
  w^\nu_\app(\xi,t) \,=\, G(\xi) + \Bigl(\frac{\nu t}{d^2}\Bigr)
  F_\nu(\xi)~, \qquad \xi \in \real^2~, \quad t > 0~,
\end{equation}
where $G$ is the Gaussian profile \eqref{Gdef} and the first
order correction $F_\nu$ is constructed as follows. Let 
$\cL$ be the Fokker-Planck operator
\[
  \cL \,=\, \Delta_\xi  + \frac12 \xi\cdot \nabla_\xi  + 1~,
  \qquad \xi \in \real^2,
\]
and $\Lambda$ be the linearized operator \eqref{Lamdef} with $w_* = G$, 
namely
\[
  \Lambda w \,=\, v^G \cdot\nabla w + K[w]\cdot\nabla G~.
\]
Here we use the functional setting of Section~\ref{sec2} in the
particular case where the asymptotic profile $w_*$ is the Oseen vortex
$G$. This means that $q(r) = \frac14 e^{-r/4}$ in \eqref{qQdef}, and
the assumptions \eqref{SS}, \eqref{qdecay}, \eqref{q2decay} are
clearly satisfied.  With these notations, the profile $F_\nu$ is the
unique solution of the linear equation
\begin{equation}\label{Fnudef}
  \frac{\nu}{\alpha}(1 - \cL) F_\nu + \Lambda F_\nu + A \,=\, 0~,
\end{equation}
where $A(\xi) = \frac{1}{2\pi}\xi_1\xi_2G(\xi)$, see \cite[Section~3.3]{Ga}. 
In polar coordinates $\xi = (r\cos\theta,r\sin\theta)$, we thus have
\[
  A(\xi) \,=\, \frac{1}{16\pi^2}\,r^2 e^{-r^2/4}\sin(2\theta) \,=\,
  \frac{1}{2\pi}\,r^2 g(r)\sin(2\theta)~, 
\] 
where $g$ is defined in \eqref{phidef}. In particular, if we use 
the angular decomposition \eqref{Xdec} of the function space 
\eqref{Xdef}, we see that $A \in X_2$, and it follows that 
$F_\nu \in X_2$ too. Now, setting $\nu = 0$ in \eqref{Fnudef}, 
we obtain the simple equation $\Lambda F_0 + A = 0$, which 
coincides with \eqref{om23} for $n = 2$. Since $X_2 \cap 
\Ker(\Lambda) = \{0\}$, we conclude that $F_0 \,=\, \omega^{(2)}$. 
It is clear from \eqref{Fnudef} that the actual profile $F_\nu$ is 
close to $F_0$ if the viscosity $\nu$ is small compared to the
circulation $\alpha$ of the vortices. As a matter of fact, it
is shown in \cite[Lemma~3.5]{Ga} that
\begin{equation}\label{Fnuapp}
  \|F_\nu - F_0\|_X \,\le\, C \frac{\nu}{\nu+\alpha}~.
\end{equation}

To formulate our main approximation result, we introduce a function
space with a weaker norm than $X$. Given any $\beta > 0$, we denote by
$Z_\beta$ the space
\[
  Z_\beta \,=\, \Bigl\{w \in L^2(\real^2)\,\Big|\, \|w\|_\beta
 < \infty\Bigr\}~, \qquad \hbox{where} \quad 
 \|w\|_\beta^2 \,=\, \int_{\real^2} |w(\xi)|^2 e^{\beta|\xi|} \d \xi~.
\]
Applying Theorem~2.5 of \cite{Ga} to the particular situation 
considered here, we obtain

\begin{proposition}\label{omapprox}
Fix $T > 0$, and let $\omega^\nu(x,t)$ be the solution of the 
rotating viscous vorticity equation \eqref{Vrotvis} with initial
data $\omega_0$. There exist positive constants $K,\beta$, 
depending only on the product $\Omega T$, such that, if 
$\omega^\nu(x,t)$ is decomposed as in \eqref{omdecomp}, then the 
vorticity profiles $w_i^\nu(\xi,t)$ satisfy
\begin{equation}\label{app1}
  \max_{i=1,2} \|w_i^\nu(\cdot,t) - w^\nu_\app(\cdot,t)\|_\beta
  \,\le\, K \Bigl(\frac{\nu t}{d^2}\Bigr)^{\!3/2}~,
\end{equation}
for all $t \in (0,T]$, where $w^\nu_\app(\xi,t)$ is given by 
\eqref{wappdef}. 
\end{proposition}

This result can be reformulated in a slightly different way, 
using the approximate solutions of Euler's equation constructed
in Section~\ref{sec2}. Indeed, if we set $\epsilon = \sqrt{\nu t}$, 
and if we remember that $w_* = G$ in the present case, we see that 
the inviscid profile $w_\epsilon$ defined in \eqref{wvasym2} satisfies
\[
  w_{\sqrt{\nu t}}(\xi) \,=\, w^\nu_\app(\xi,t) + \Bigl(\frac{\nu t}
  {d^2}\Bigr)\Bigl(F_0(\xi) - F_\nu(\xi)\Bigr) + \cO\Bigl(\Bigl(
  \frac{\nu t}{d^2}\Bigr)^{\!3/2}\Bigr)~.
\]
Thus, combining \eqref{Fnuapp} and \eqref{app1}, we obtain\:

\begin{corollary}\label{cor2}
Under the assumptions of Proposition~\ref{omapprox}, we have
\begin{equation}\label{app2}
  \max_{i=1,2} \|w_i^\nu(\cdot,t) - w_{\sqrt{\nu t}}\|_\beta
  \,\le\, K \Bigl(\frac{\nu t}{d^2}\Bigr)^{\!3/2} + C\,\frac{\nu t}
  {d^2}\,\frac{\nu}{\nu+\alpha}~, 
\end{equation}
for all $t \in (0,T]$, where $w_\epsilon$ is the inviscid profile
defined by \eqref{wvasym2} with $w_* = G$.  
\end{corollary}

Finally, using the continuous inclusion $Z_\beta \hookrightarrow 
L^1(\real^2)$, we can formulate an approximation result for 
the original function $\omega^\nu(x,t)$. If we compare \eqref{Vpair2}, 
\eqref{omdecomp} and use the fact that $G(\xi)$ and $F_0(\xi)$ are 
even functions of $\xi$, we arrive at\:

\begin{corollary}\label{cor3}
Under the assumptions of Proposition~\ref{omapprox}, we have
\begin{equation}\label{app3}
  \frac{1}{\alpha}\int_{\real^2} | \omega^\nu(x,t) - 
  \omega_{\sqrt{\nu t}}(x)|\d x \,\le\, K \Bigl(\frac{\nu t}{d^2}
  \Bigr)^{3/2} + C\,\frac{\nu t}{d^2}\,\frac{\nu}{\nu+\alpha}~, 
\end{equation}
for all $t \in (0,T]$, where $\omega_\epsilon$ is the approximate
steady state of the rotating Euler equation defined in \eqref{Vpair2}
with $w_* = G$.
\end{corollary}

Comparing this last result with Theorem~\ref{thm1}, we see that
replacing a linear superposition of Oseen vortices by a more accurate
solution of Euler's equation, which takes into account the deformation
of the vortex cores due to mutual interaction, results in a better
approximation. If we believe that there exists an {\em exact} solution 
of Euler's equation which is close to $\omega_\epsilon$ for $\epsilon$
sufficiently small, then \eqref{app3} shows that the solution
$\omega^\nu(x,t)$ of the rotating viscous vorticity equation
\eqref{Vrotvis} slowly travels through a family of steady states of 
the inviscid equation \eqref{Vrot} indexed by the length
parameter $\epsilon$, which evolves diffusively according to $\epsilon
= \sqrt{\nu t}$.

\begin{remark}\label{final}
The right-hand side of \eqref{app3} suggests that both dimensionless
quantities $\nu/\alpha$ and $\nu t/d^2$ play an important role 
in the evolution of viscous vortices. It is tempting to eliminate
one of these quantities by considering, for instance, the limit 
$\nu \to 0$ while $\epsilon = \sqrt{\nu t}$ is kept fixed. 
Under the assumptions of Proposition~\ref{omapprox}, one may
conjecture that the solution $\omega^\nu(x,t)$ of \eqref{Vrotvis}
satisfies, if $\epsilon > 0$ is sufficiently small, 
\[
  \omega^\nu(x,\epsilon^2/\nu)  ~\xrightarrow[\nu \to 0]{}~ 
  \omega_\epsilon(x)~, \qquad \hbox{for all~} 
  x \in \real^2~,
\]
where $\omega_\epsilon(x)$ is an exact stationary solution of 
\eqref{Vrot} of the form \eqref{Vpair2}. Unfortunately, the 
results of \cite{Ga} do not provide any control on $\omega^\nu(x,t)$ 
in the limit where $\Omega t \to \infty$. 
\end{remark}

\end{document}